\magnification 1200
        \def\R{{\rm I\kern-0.2em R\kern0.2em \kern-0.2em}}
        \def\N{{\rm I\kern-0.2em N\kern0.2em \kern-0.2em}}
        \def\P{{\rm I\kern-0.2em P\kern0.2em \kern-0.2em}}
        \def\B{{\rm I\kern-0.2em B\kern0.2em \kern-0.2em}}
        \def\Z{{\rm I\kern-0.2em Z\kern0.2em \kern-0.2em}}
        \def\C{{\bf \rm C}\kern-.4em {\vrule height1.4ex width.08em depth-.04ex}\;}
        \def\B{{\bf \rm B}\kern-.4em {\vrule height1.4ex width.08em depth-.04ex}\;}
        
        \def\D{{\Delta}}

        \def\z{{\zeta}}

        \
        \vskip 10mm
        \centerline {\bf ON DISCS IN BIDISCS}
        \vskip 4mm
        \centerline{Josip Globevnik}
        \vskip 4mm

 \noindent        \bf Abstract\ \ \rm Let $\Delta$ be the open unit disc in $\C $. We show that there is no continuous map $F\colon\overline \Delta\rightarrow \overline\Delta^2$, holomorphic on $\Delta$ and such that $F(b\Delta) =  b(\Delta^2). $
        \vskip 6mm
 Denote by $\D$    the open unit disc in $\C$. Given a bounded convex domain $D\subset \C^2$ we consider holomorphic maps $F\colon\D\rightarrow D$ which extend continuously to $\overline\D$ and which are proper, that is, satisfy $F(b\D )\subset bD. $   We ask how large $F(b\D)$ can be.  In particular, can we have $F(b\D)= bD?$ It is known that the answer is positive in the case of the ball $D=\B = \{ (z,w)\colon |z|^2+|w|^2<1\}$:
\vskip 2mm
\noindent\bf PROPOSITION 1\ \rm\ [G, Cor.\ 2] \it There is a continuous map $F\colon\overline\D\rightarrow \overline \B$, holomorphic on $\D$, and such that $F(b\D)=b\B$.\rm 
\vskip 2mm
\noindent In the present note we show that in general the answer to the preceding  question is no. In particular, it is negative for $D$ = $\D^2$:
\vskip 2mm
\noindent \bf PROPOSITION 2\ \it \ There is no continuous map $F\colon\overline\D\rightarrow \overline\D^2 $, holomorphic on $\Delta$ and such that $F(b\D)= b(\D^2)$.
\vskip 2mm
\noindent\bf DEFINITION\ \ \it A set of the form $\{\z\}\times \D$ or $\D\times\{\z\}$ where $\z\in b\D$ is called \rm an open face \it of the bidisc $\D^2$.  \rm
\vskip 2mm
\noindent \bf Proof of Proposition 2. \rm Let $F=(f,g)\colon\overline \D\rightarrow \overline\D^2$ be a continuous map, holomorphic on $\D$ and such that $F(b\D)\subset b(\D^2).$ 

If one of the components, say $f$, is a constant $\alpha$ then $F(b\D)\subset \{ \alpha\}\times \overline\D$  so $F(b\D)$ cannot equal $b(\Delta^2)$.  So assume that both $f,g$ are nonconstant. 
We shall show that for any open face $\Phi$ of $\D^2$ the set $F(b\D )\cap\Phi$ has no cluster point in $\Phi$.  In fact, we show that

$$
\left .\eqalign
{\hbox{\ given\ } \xi\in b\D \hbox{\ there is no injective sequence\ } \z_n\in b\D \hbox{\ such that\ }  f(\z_n)=\xi \cr \hbox{\ for all\ }n  
\hbox{\ and such that the sequence\ } g(\z_n) \hbox {\ has a cluster point in\ } \D. \cr }\right \}. \ \ \eqno (1)
$$
\rm Assume the contrary, so that $\z_n\in b\D$ is an injective sequence such that $f(\z_n ) = \xi$ for all $n$ and such that $g(\z _n)$ has a cluster point in $\D$. Passing to a subsequence we may assume that $ \z_n$ converges to $\z_0\in b\D$ and that $\lim g(\z_n) = g(\z_0)\in\D$.  Passing to a subsequence we may assume that all $\z_n$ are contained in a small open arc $J\subset b\D$ centered at $\z_0$ such that $|g(\z )|<1\ \ (\z\in J)$. Since $F(b\D)\subset b(\D^2)$ it follows that $|f(\z)|=1\ \ (\z\in J)$. Denote by $^\ast$ the reflection across $b\D$ so for $z\in\D\setminus\{ 0\} $ write $z^\ast = 1/\overline z$. By the reflection principle there is a narrow open  neighbourhood $U$ of $J$ in $\C$ such that $U^\ast = U$ and such that by defining $f$ on $U\cap[\C\setminus \overline \D]$ by
$$
f(w) = [f(w^\ast)]^\ast \ \ \ \ (w\in U\cap (\C\setminus \overline \D) 
$$ 
 the function $f $ is holomorphic on $\D\cup U$ and satisfies
$$
|f(\z)|<1 \ \ (\z\in U\cap \D) \hbox{\ \ and\ \ } |f(\z)|>1  \  \ (\z\in U\cap[\C\setminus\overline\D] ).  
\eqno (2)
$$
We show  that the derivative $f^\prime $ has no zero on $J$.  Indeed, if $f^\prime (z)=0$ for some $z\in J$ then since $ f$ is not a constant, $z$ is an isolated zero of $f^\prime$ and the function $w\mapsto f(w)-f(z)$  has a zero of order at least two at $z$. This means that for every sufficiently small circle $\gamma$ centered at $z$ and contained in $U$ the winding number of $f\circ \gamma$ around $f(z)$ is at least two which is impossible since by (2) $f$ maps $\gamma\cap\D $ to $\D$ and $\gamma\cap [\C\setminus\D]$ to $\C\setminus \D$. 
This implies that $f$ is locally one to one on $J$. In particular, $f$ is one to one in a neighbourhood of $\z_0$ in $J$  which contradicts the fact that $f(\z _n)=\xi$ for all $n$ and that the injective sequence $\z_n$ converges to $\z_0$. 
  This proves (1) and the same holds with the roles of $f$ and $g$ interchanged. This completes the proof of Proposition 2.
\vskip 3mm
\noindent\bf REMARK  \ \rm If one drops the requirement about boundary continuity one can do much more. For instance, one can show that given any bounded convex domain $D$ there is a proper holomorphic embedding $F\colon\D\rightarrow D$ such that $\overline{F(\D)} \setminus F(\D)=bD$  \ [FGS] .
\vskip 10mm

\centerline {\bf REFERENCES}
\vskip 4mm

\noindent [FGS] F.\ Forstneri\v c, J.\ Globevnik, B. Stensones: Embedding holomorphic discs through discrete sets. Math.\ Ann,\ 305 (1996)559-569
\vskip 2mm 
\noindent [G]  Globevnik, J.: Boundary interpolation and proper holomorphic maps from the disc to the ball.  Math\ Z.\ 198 (1988) 143-150
\vskip 8mm
\noindent Institute of Mathematics, Physics and Mechanics and

\noindent  Department of Mathematics, University of Ljubljana

\noindent Jadranska 19, 1000 Ljubljana, Slovenia  

\noindent josip.globevnik@fmf.uni-lj.si

\

\bye